\documentclass[a4paper,12pt,leqno]{article}
\usepackage{epic,eepic}
\usepackage{macros}
\usepackage{latexsym}
\begin{document}
\title{ 
  Metrics of constant curvature $1$ with\\
  three conical singularities on $2$-sphere 
}
\date{January 28, 1998}
\author{\normalsize \bf
   \begin{tabular}[t]{c}
    Masaaki Umehara\\
    {\footnotesize\rm Osaka Univ.}
   \end{tabular}
   and 
   \begin{tabular}[t]{c}
    Kotaro Yamada%
    \raisebox{0.8ex}{\tiny\rm 1}\\
    {\footnotesize\rm Kumamoto Univ.}
   \end{tabular}
}
\maketitle
\renewcommand{\thefootnote}{\relax}
\addtocounter{footnote}{1}
\renewcommand{\thefootnote}{\arabic{footnote}}
\footnotetext{Supported by Sumitomo Foundation.}

\section*{Introduction.}

Conformal metrics of constant curvature $1$ with conical 
singularities on a closed Riemann surface $\Sigma$
is bijectively corresponding to branched \cmcone{} surfaces in 
the hyperbolic 3-space with given hyperbolic Gauss map
defined on ${\Sigma}$ excluded finite points (see \cite{UY3} 
and also \cite{ruy2}).

It is a fundamental problem 
to  consider the existence of such metrics.
Several important results for this direction
are found in \cite{Troyanov1}, \cite{CL}, \cite{LT}.
Troyanov \cite{Troyanov2} gave a classification of metrics
of constant curvature $1$ with at most two conical singularities
on $2$-sphere.

In this paper, we shall give a necessary and sufficient condition
for the existence and the uniqueness of metric with 
three conical singularities 
of given order
on 2-sphere.
In Section \ref{sec:preliminaries}, we recall some basic properties of
null meromorphic curves in $\PSL(2,\C)$. 
The irreducible case is proved in Section \ref{sec:irreducible}  
using the study of \cmcone{} surfaces developed in \cite{ruy1}.
Moreover, applying the same method, 
we classify all genus zero 
irreducible \cmcone{} surfaces with embedded three regular ends.
In Section \ref{sec:reducible}, we give 
a method for explicit construction of
all reducible metric with three singularities.

The authors thank Wayne Rossman and Mikio Furuta for
fruitful discussions and encouragement.
They are also wishing to thank Yoshinobu Hattori for informative 
conversations.

\section{Preliminaries.}
\label{sec:preliminaries}

In this section, we recall fundamental properties of
null meromorphic curves in $\PSL(2,\C)$.

\begin{Def}
  Let $F: \Sigma \to \PSL(2,\C)$ be a meromorphic map defined
  on Riemann surface $\Sigma$.
  Then $F$ is called {\it null\/} if
  \[
            \det(F^{-1}\cdot F_z)=0 
  \]
  holds on $\Sigma$, where $z$ is a complex coordinate.
  (The condition does not depend on the choice of coordinates.) 
\end{Def}

Let $F: \Sigma \to \PSL(2,\C)$ be a null meromorphic map. 
We define a matrix $\alpha$ by
\[
   \alpha=\pmatrix{ \alpha_{11} & \alpha_{12} \cr
                        \alpha_{21} & \alpha_{22} }
         :=     F^{-1}\cdot dF,
\]
and set
\begin{equation}\label{eq:gomega}
   g:= \alpha_{11}/\alpha_{21},\qquad
  \omega:=\alpha_{21}.
\end{equation}
Then the pair $(g,\omega)$ of a meromorphic function $g$
and a meromorphic $1$-form $\omega$ on $\Sigma$ satisfies the equality
\begin{equation}\label{eq:ode0}
   F^{-1}\cdot dF= \pmatrix{ 
                        g & -g^2 \cr
                        1 & -g\hphantom{^2}} \, \omega.
\end{equation}
Conversely, let $g$ be a meromorphic function and $\omega$
a meromorphic $1$-form on $\Sigma$. Then the ordinary
differential equation \eqref{eq:ode0} is integrable and the solution
$F$ is a null map into $\PSL(2,\C)$, however
$F$ may not be single-valued on $\Sigma$. Moreover,
$F$ may have essential singularities.
We call the pair $(g,\omega)$ the {\it Weierstrass data\/} of $F$. 

\begin{Def}
 Let 
 \[
     F=\pmatrix{F_{11} & F_{12} \cr F_{21} & F_{22}}
 \]
 be a null meromorphic map of $\Sigma$ into $\PSL(2,\C)$.
 We call 
 \[
      G=\frac {dF_{11}}{dF_{21}}=\frac {dF_{12}}{dF_{22}}
 \]
 the {\it hyperbolic Gauss map\/} of $F$.
 Furthermore, we  call $g$ in \eqref{eq:gomega} the {\it secondary Gauss
 map\/}
 and $Q=\omega dg$ the {\it Hopf differential\/} of $F$.
\end{Def}

Let $F: \Sigma \to \PSL(2,\C)$
be a null meromorphic map.
Then for $a,b\in \SL(2,\C)$,
$\hat F =a\cdot F\cdot b^{-1}$ is also a null meromorphic map.
The associated two Gauss maps $\hat G$, $\hat g$, and
the Hopf differential $\hat Q$ of $\hat F$ are given by
\begin{equation}\label{eq:three}
  \hat G=\frac{a_{11}G+a_{12}}{a_{21}G+a_{22}}, \quad
  \hat g=\frac{b_{11}g+b_{12}}{b_{21}g+b_{22}}, \quad
  \mbox{and}\quad
  \hat Q=Q, 
\end{equation}
where we set
\[
  a=\pmatrix{ a_{11}& a_{12}\cr
              a_{21}& a_{22}},
  \quad b=\pmatrix{ b_{11}& b_{12}\cr
           b_{21}& b_{22}}.
\]

Let $(U;z)$ be a complex coordinate of $\Sigma$.
Now we consider the Schwarzian
derivatives
 $S(G)$ and $S(g)$ of $G$ and $g$ respectively, 
where
\begin{equation}\label{eq:sch}
  S(G)=
  \left[
  \left(\frac{G''}{G'}\right)'-\frac 12 \left(\frac{G''}{G'}\right)^2
  \right]\,dz^2 \qquad
  \left( '=\frac{d}{dz}\right).
\end{equation}
The description of the Schwarzian derivative depends on
the choice of complex coordinates.
However, any difference of two Schwarzian derivatives
does not depend on the choice of
complex coordinate. The following identity can be checked:
\begin{equation}\label{S-rel}
  S(G)= S(g) - 2Q. 
\end{equation}

Conversely, the following fact is known:

\begin{Fact}[\cite{Sm}, \cite{UY3}]\label{fact1}
 Let $G$ and $g$ {\rm ($S(G)\not\equiv S(g)$)} be non-constant
 meromorphic functions
 on a Riemann surface $\Sigma$. 
 Then there exists a unique
 null meromorphic map $F:\Sigma\to \PSL(2,\C)$ 
 with the hyperbolic Gauss map $G$ and the Hopf differential 
 \[
    Q=-\frac 12(S(G)-S(g)),
 \]
 such that $(g,Q/dg)$ is the Weierstrass data of $F$.
\end{Fact}

Moreover, the following fact plays an important role in the latter 
discussions:

\begin{Fact}[\cite{UY4}]\label{fact2}
 Let $F:\Sigma\to \PSL(2,\C)$ be 
 a null meromorphic map with hyperbolic Gauss map $G$
 and secondary  Gauss map $g$.
 Then the inverse map $F^{-1}$ is a 
 null meromorphic map with hyperbolic Gauss map $g$
 and secondary  Gauss map $G$. In particular, $F$ 
 satisfies the following ordinary differential equation
 \[
     -F\cdot d(F^{-1})=dF\cdot F^{-1}=
     \pmatrix{ G & -G^2 \cr
               1 & -G\hphantom{^2}  } \, \frac {Q}{dG}.
 \]
\end{Fact}

Finally, we point out the following elementary fact on linear
algebra:

\begin{Fact}[\cite{ruy1}]\label{fact3}
  A matrix $a\in \SL(2,\C)$ satisfies $a\cdot \bar a=\id$
  if and only if $a$ is of the form
  \[
      a=\pmatrix{       p & i \gamma_1 \cr
                        i\gamma_2 & \bar p },
      \qquad (\gamma_1,\gamma_2\in \R, \quad p\bar p+\gamma_1\gamma_2=1).
  \]
Moreover, $a$ can be 
diagonalized by a real matrix in $\SL(2,\R)$
whenever it is semi\-simple.
\end{Fact}
\section{Irreducible metrics with three singularities.}
\label{sec:irreducible}

A conformal positive semidefinite metric $d\sigma^2$ on a closed 
Riemann surface $\Sigma$ is said to have a {\it conical singularity\/} of 
order $\beta$ ($>-1$) at $p$ on $\Sigma$ 
if it is asymptotic to the metric of 
the form $C |z-p|^{2\beta}\,dz\,d\bar z$ at $p$, 
where $C$ is a constant and $z$ is a complex local coordinate
of $\Sigma$. 
We denote by $\Met(\Sigma)$ the set of (non-vanishing) 
conformal pseudometrics of constant curvature $1$ 
on $\Sigma$ with finitely many conical singularities.
Let $p_1$, \ldots, $p_n$ be the set of singular points of the given
pseudometric $d\sigma^2\in \Met(\Sigma)$ and $\beta_j$
the order of conical singularity at $z=p_j$.
We call a formal sum
\[
    \beta_1 p_1+\cdots +\beta_n p_n
\]
the {\it divisor\/} of $d\sigma^2$.

For each metric $d \sigma^2\in \Met(\Sigma)$, there exists
a meromorphic function $g$ defined on the universal cover 
$\tilde\Sigma_{p_1,\ldots,p_n}$ of 
$\Sigma_{p_1,\ldots,p_n}:=\Sigma\setminus \{p_1,\ldots,p_n\}$ 
such that the metric is the pull-back of the canonical metric
$d\sigma^2_0:=4\,dz\,d\bar z/(1+|z|^2)^2$ on $S^2=\C\cup\{\infty\}$,
that is, we have an expression
\begin{equation}\label{eq:pullback}
   d\sigma^2=g^*d\sigma^2_0=\frac{4\,dg\,d\bar g}{(1+|g|^2)^2}.
\end{equation}
Such a function $g$ is unique up to the change
\begin{equation}\label{eq:g-change}
    g \mapsto a\star g, \qquad 
    \mbox{where}
    \quad
    a\in\PSU(2):=\SU(2)/\{\pm\id\}.
\end{equation}
Here, for $a\in\PSL(2,\C)$ in general, we set
\begin{equation}\label{eq:moebius}
      a\star g =
      \frac{a_{11}g+a_{12}}{a_{21}g+a_{22}},
      \qquad 
      \mbox{where}
      \quad
         a=\pm\pmatrix{
                      a_{11} & a_{12} \cr
                      a_{21} & a_{22}}.
\end{equation}

We denote by  
$\pi:\tilde \Sigma_{p_1,\ldots,p_n}\to \Sigma_{p_1,\ldots,p_n}$ 
the covering projection. 
Fix a base point $\tilde z_0$ on 
$\tilde \Sigma_{p_1,\ldots,p_n}$ and set $z_0=\pi(\tilde z_0)$.
For each $z\in \pi^{-1}(z_0)$, there exists 
a unique deck transformation $T$ such that $T(\tilde z_0)=z$.
Thus the fundamental group $\pi_1(\Sigma_{p_1,\ldots,p_n})$ is 
identified with
the deck transformation group.
By \eqref{eq:pullback} and \eqref{eq:g-change}, 
there exists a representation 
$\rho_g :\pi_1(\Sigma_{p_1,\ldots,p_n})\to \PSU(2)$
such that 
\begin{equation}\label{eq:g-repr}
   g\circ T^{-1}=\rho_g(T)\star g,
   \qquad
   \left(T \in\pi_1(\Sigma_{p_1,\ldots,p_n}^{})\right).
\end{equation}
Later, we will see that the representation $\rho_g$ can be
lifted to an
$\SU(2)$-re\-pre\-sen\-ta\-tion.

Metrics $d\sigma^2 \in \Met(\Sigma)$
are divided into the three classes defined below;
\begin{enumerate}
  \item\label{item:red:1}
        A metric is called {\it irreducible\/} when the image of 
        the representation $\rho_g$ is not abelian.
  \item\label{item:red:2} 
        A metric is called {\it $\Hyp^1$-reducible\/} when the
        image of the representation $\rho_g$ is abelian but non-trivial.
  \item\label{item:red:3} 
        Then a metric is called {\it $\Hyp^3$-reducible\/} when
        the image of the representation $\rho_g$ is trivial.
\end{enumerate}
\label{page:reducible}
If there exists $a\in \PSL(2,\C)$ such that the 
image of the representation $a\cdot \rho_g\cdot a^{-1}$ 
is also contained in $\PSU(2)$,
then another metric $d\sigma^2_a:=(a\star g)^*d\sigma^2_0$ has the same 
divisor 
and the  Schwarzian derivative as $d\sigma^2$.
Hence, the set
\begin{equation}\label{eq:deform-sp}
  \left\{
    d\sigma^2_a=(a\star g)^*d\sigma_0^2
    \,|\,a\in\PSL(2,\C)\,;\,a\cdot\Im\rho_g\cdot a^{-1}\subset \PSU(2)
  \right\}
\end{equation}
defines a deformation space of the metric $d\sigma^2=d\sigma^2_{\id}$
preserving the divisor
and the Schwarzian derivative.
Since $d\sigma_a^2=d\sigma^2$ for $a\in\PSU(2)$,
the set \eqref{eq:deform-sp} is identified with the set $I_{\Gamma}$ 
in Appendix~\ref{app:B}, where $\Gamma=\Im\rho_g$.

Suppose that $d\sigma^2$ is $\Hyp^1$-reducible
(resp.\ $\Hyp^3$-reducible). 
Then  by Lemma~\ref{lem:B} in Appendix~\ref{app:B}, 
the set \eqref{eq:deform-sp} is identified 
with a totally geodesic subset 
of dimension one (resp.\ three) in the hyperbolic $3$-space.

\medskip

We now determine all of irreducible metrics in $\Met(S^2)$
with three conical singularities.
Reducible case is discussed in the next section. 

\begin{Thm}\label{thm:irred}
 There exists an irreducible metric $d\sigma^2$ in $\Met(S^2)$
 with a given divisor
 \begin{equation}\label{eq:divisor}\label{eq:D}
    D:=\beta_1 p_1 +\beta_2 p_2 + \beta_3 p_3
       \qquad (\beta_j>-1)
 \end{equation}
 if and only if the following inequality
 holds{\rm :}
 \begin{equation}\label{eq:condition}\label{eq2}
    \cos^2 B_1 +\cos^2 B_2 + \cos^2 B_3+
               2\cos B_1\cos B_2 \cos B_3<1,
 \end{equation}
 where $B_j:=\pi(\beta_j+1)$ $(j=1,2,3)$.
 Moreover, such a metric $d\sigma^2$ is uniquely determined.
\end{Thm}

\begin{Remk}
 The condition \eqref{eq:condition} implies that all $\beta_j$ 
 ($j=1,2,3$) are not integer.
 In fact, suppose one of $\beta_j$'s, say $\beta_1$, is an integer.
 Then $\cos B_1=\pm 1$, and hence \eqref{eq:condition}
 fails. Conversely, if 
all $\beta_j$ 
 ($j=1,2,3$) are not integer, the metric is automatically irreducible
(see Corollary \ref{added}).
\end{Remk}

\begin{Remk}
 Recently, Hattori \cite{H} gave an alternative approach to
 prove Theorem~\ref{thm:irred} based on the geometry of spherical 
 polytopes.
\end{Remk}

\begin{Remk}

 If a metric $d\sigma^2\in \Met(\Sigma)$ with divisor $D$ in \eqref{eq:D}
 is reducible, then Lemma~\ref{lem:A} in Appendix~\ref{app:A}
 yields the equality
 \begin{equation}\label{eq:eq}
    \cos^2 B_1 +\cos^2 B_2 + \cos^2 B_3+
               2\cos B_1\cos B_2 \cos B_3=1.
 \end{equation}
 However, even if $D$ satisfies \eqref{eq:eq},
 it does not imply the existence of the metric.
 In fact, if $\beta_1,\beta_2,\beta_3\not\in \Z$,
 such a metric never exists (see Lemma \ref{lem:integer} and Corollary~\ref{added}).
\end{Remk}

We identify $S^2$ with $\C\cup \{\infty\}$ by
the stereographic projection, and let $z$ 
be the canonical complex coordinate of $\C$.
For each $d\sigma^2\in \Met(S^2)$,  
we define the Schwarzian derivative $\tilde S(d\sigma^2)$ as
\begin{equation}\label{eq:schwarz}\label{eq3}
  \tilde S(d\sigma^2):=S(g)=
  \left[
  \left(\frac{g''}{g'}\right)'-\frac 12 \left(\frac{g''}{g'}\right)^2
  \right]\,dz^2 \qquad
  \left( '=\frac{d}{dz}\right),
\end{equation}
where $g$ is a meromorphic function defined on  
$\tilde S^2_{p_1,p_2,p_3}$ satisfying \eqref{eq:pullback}.
This is independent of the choice of $g$.
In fact $g$ is unique up to change as in \eqref{eq:g-change},
and a M\"obius transformation \eqref{eq:moebius} on $g$ does not change 
$S(g)$.

Let $d\sigma^2\in\Met(S^2)$ with divisor $D$ as in \eqref{eq:D},
and take a function $g$ as in \eqref{eq:pullback}.
Then for each  $p_j$,
 there exists $a\in\PSU(2)$ such that
\begin{equation}\label{eq:g-expand}
  a\star g = (z-p_j)^{\beta_j+1}(g_0 + g_1 z + \cdots ),\qquad (g_0\neq 0).
\end{equation}
Hence $\tilde S(d\sigma^2)$ can be written in the following leading
terms in the Laurent expansions at $z=p_j$:
\begin{equation}\label{eq:schwarz-exp}
  \tilde S(d\sigma^2)= 
   \left[
     -\frac{\beta_j(\beta_j+2)}{2}\frac{1}{(z-p_j)^2}+\cdots
   \right]\,dz^2,
\end{equation}
and $\tilde S(d\sigma^2)$ is holomorphic on $S^2_{p_1,p_2,p_3}$.
Since the M\"obius transformation group acts on the sphere, 
we may assume 
\begin{equation}\label{eq:pos-sing}
     p_1=0,\qquad
     p_2=1,\qquad
     \mbox{and}\qquad
     p_3=\infty
\end{equation}
without loss of generality.
By \eqref{eq:schwarz-exp}, 
$\tilde S(d\sigma^2)$ is uniquely determined,
since the total order of a holomorphic 2-differential 
on $S^2$ is $4$.
In fact,
We have
\begin{equation}\label{eq:schwarz-2}
  \tilde S(d\sigma^2)=\left[
               \frac{c_3 z^2+(c_2-c_1-c_3)z+c_1}{z^2(z-1)^2}
            \right]\,dz^2,
\end{equation}
where $c_j=-\beta_j(\beta_j+2)/2\in\R$  ($j=1, 2, 3$).
Now we set
\begin{equation}\label{eq:GQ}
  G:=z, \qquad
  \mbox{and}
  \qquad
  Q:=
   \frac{1}{2}
    \left(\frac{c_3z^2+(c_2-c_1-c_3)z+c_1}{z^2(z-1)^2}\right)
     dz^2.
\end{equation}
By Fact \ref{fact1},
there exists a unique null holomorphic map $F_g:
\tilde S^2_{p_1,p_2,p_3}\to \PSL(2,\C)$
such that $G\circ \pi$ and $g$ are the hyperbolic Gauss map
and the secondary Gauss map respectively.
Since $\tilde S^2_{p_1,p_2,p_3}$ is simply connected,
$F_g$ can be lifted to a null holomorphic map
$\tilde F_g:
\tilde S^2_{p_1,p_2,p_3}\to \SL(2,\C)$.
By Fact \ref{fact2},
we have the following relations:
\begin{eqnarray}
\label{eq:ode}
  & & d\tilde F_g
         \cdot \tilde F_g^{-1}= 
         \pmatrix{G &-G^2 \cr
                  1 & -G\hphantom{^2}} \frac{Q}{dG}, \\
\label{eq:g-def}
  & &  g = -\frac{d\tilde F_{12}}{d\tilde F_{11}}=
           -\frac{d\tilde F_{22}}{d\tilde F_{21}},
\end{eqnarray}
where $\tilde F_g=(\tilde F_{ij})$.
The  right-hand side of \eqref{eq:ode} is single-valued
and has poles at $p_1$, $p_2$ and $p_3$.
Thus, there exists a representation 
$\rho_{\tilde F_g}:\pi_1(S^2_{p_1,p_2,p_3})\to  \SL(2,\C)$
such that 
\begin{equation}\label{eq:repr}
  \tilde F_g\circ T=\tilde F_g\cdot \rho_{\tilde F_g}(T) \qquad 
   \left(T\in 
        \pi_1(S^2_{p_1,p_2,p_3})
  \right).
\end{equation}
By \eqref{eq:three}, we have 
\begin{equation}\label{eq:rho-g-rho-F}
  \rho_g(T)=\pm\rho_{\tilde F_g}(T)\in\PSU(2)
  \qquad \left(
            T\in\pi_1(S^2_{p_1,p_2,p_3})
         \right).
\end{equation}
Here, we consider $\pm a$ ($a\in\SU(2)$) as an element of
$\PSU(2)=\SU(2)/\{\pm\id\}$.
In particular, $\rho_{\tilde F_g}$ is an 
$\SU(2)$-representation 
of $\pi_1(S^2_{p_1,p_2,p_3})$.

Conversely, if $\tilde F=(\tilde F_{ij})$ is a solution of 
\eqref{eq:ode} such that the image of $\rho_{\tilde F}$ lies in $\SU(2)$, then
\eqref{eq:rho-g-rho-F} holds for a meromorphic function $g$
defined by \eqref{eq:g-def}. 
So $d\sigma^2$ defined by \eqref{eq:pullback} is well-defined on $S^2$,
that is $d\sigma^2\in\Met(S^2)$, which has divisor $D$
in \eqref{eq:D}.
Let $T_j$ be the deck transformation on $\tilde S^2_{p_1,p_2,p_3}$
induced from the loop at $p_j$ ($j=1, 2, 3$).
By \eqref{eq:g-expand} and \eqref{eq:rho-g-rho-F}, the eigenvalues
of $\rho_{\tilde F_g}(T_j)$ are 
$\{\pm e^{iB_j}, \pm e^{-iB_j}\}$, where $B_j=\pi(\beta_j+1)$.
Moreover, we can eliminate this $\pm$-ambiguity as follows:

\begin{Lemma}\label{lem:eigen}
  The eigenvalues of $\rho_{\tilde F_g}(T_j)$ are $\{-e^{iB_j}, -e^{-iB_j}\}$,
  where $B_j=\pi(\beta_j+1)$ $(j=1, 2, 3)$.
  In particular, $\tr \rho_F(T_j)=-2\cos B_j$ holds.
\end{Lemma}
\begin{Proof}
  Let $\tilde F_{D}$ be the solution of ordinary differential equation
  \[
      \tilde F_D \cdot d\tilde F_D^{-1}=-d\tilde F_D
        \cdot \tilde F_D^{-1}= 
     -\pmatrix{G &-G^2 \cr
               1 & -G\hphantom{^2}} \frac{Q}{dG},
      \qquad \tilde F_D(\tilde z_0)=\id
  \]
  for the divisor $D$, where $\tilde z_0\in \tilde S^2_{p_1,p_2,p_3}$ 
  is a base point.
  Then the monodromy representation
  $\rho_D:\pi_1(S^2_{p_1,p_2,p_3})\to \SL(2,\C)$ with respect to $\tilde F_D$
  is conjugate to the representation $\rho_{\tilde F_g}$.
  Let
  $\tau_j:=\tr \rho_D(T_j)=\pm 2\cos B_j=\pm 2 \cos\pi(\beta_j+1)$.
  Since \eqref{eq:ode} is real analytic in parameters
  $(\beta_1,\beta_2,\beta_3)$, 
  $\tau_j$'s 
  are  real analytic functions 
  in $\beta_j$.
  Here, let $\beta_1=\beta_2=\beta_3=0$.
  In this case,  $\tilde F_D$ is constant because $Q\equiv 0$,
  and hence each 
  $\rho_D(T_j)$ is the identity matrix.
  Thus we have
  \[
     \tau_j|_{(\beta_1,\beta_2,\beta_3)=(0,0,0)}=2=
             -2\cos B_j|_{(\beta_1,\beta_2,\beta_3)=(0,0,0)},
  \]
  and, by real analyticity, $\tau_j=-2\cos B_j$.
  Hence $\tr \rho_{\tilde F_g}(T_j)=\tr \rho_{D}(T_j) = -2 \cos B_j$, and
  the eigenvalues of $\rho_{D}(T_j)$ are $-e^{\pm iB_j} $.
\end{Proof}

\begin{Proof*}{Proof of Theorem\/ {\rm \ref{thm:irred}}}
 First, we prove the uniqueness  of the irreducible metric as follows:
 Let $d\sigma_l^2\in \Met(S^2)$
 ($l=1,2$) be  
 two irreducible metrics with
 the same divisor $D$ as in \eqref{eq:D}.
 Then we have 
 \begin{equation}\label{eq:sigma-eq}
     2Q=\tilde S(d\sigma^2_1)=\tilde S(d\sigma^2_2).
 \end{equation}
 Let $g_l$  be 
 meromorphic functions on $\tilde S^2_{p_1,p_2,p_3}$
 such that $d\sigma^2_l=g_l^*d\sigma^2_0$ ($l=1,2$).
 By \eqref{eq:ode}, $\tilde F_{g_1}$ and
 $\tilde F_{g_2}:\tilde S^2_{p_1,p_2,p_3}\to \SL(2,\C)$ 
 satisfy the same ordinary differential equation.
 Thus,
 there exists a matrix $a\in\SL(2,\C)$ such that 
 $\tilde F_{g_2}=\tilde F_{g_1}\cdot a$.
 So we have
\begin{equation}\label{eq:rho1-rho2}
  \rho_{\tilde F_{g_2}}^{}=a^{-1}\cdot \rho_{\tilde F_{g_1}}^{}\cdot a
  \qquad 
  \left(a\in\SL(2,\C)\right).
\end{equation}
Since $\rho_{\tilde F_{g_1}}$ and $\rho_{\tilde F_{g_2}}$ are 
both $\SU(2)$-representations and they are irreducible, 
we have $a\in\SU(2)$ (see Appendix~\ref{app:B}).
This implies $d\sigma_1^2=d\sigma_2^2$.

\medskip

Next, we take $d\sigma^2\in\Met(S^2)$ with divisor $D$ in \eqref{eq:D}.
Let $T_j$ ($j=1, 2, 3$) be a deck transformation
corresponding to the loop at $z=p_j$.
By Lemma \ref{lem:eigen}, 
the matrix $\rho(T_j)$ ($j=1, 2, 3$) must have
eigenvalues of the form $- e^{\pm i B_j}$.
Since $\rho(T_1)\cdot \rho(T_2)\cdot \rho(T_3)=\id$, 
we have the inequality
\[
\cos^2 B_1 +\cos^2 B_2 + \cos^2 B_3+
               2\cos B_1\cos B_2 \cos B_3\leq 1
\]
by Lemma~\ref{lem:A} in Appendix~\ref{app:A}. 
Since $d\sigma^2$ is irreducible, equality never holds.
 
\medskip

Finally,
we show the existence of the metric under the condition 
\eqref{eq:condition}.
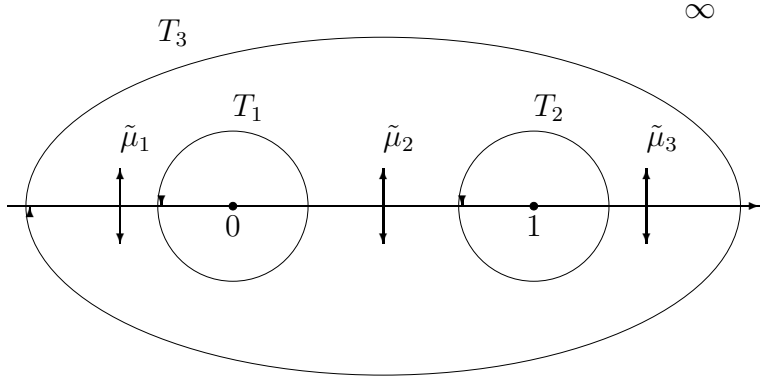
\begin{figure}
\begin{center}
\setlength{\unitlength}{1cm}
\begin{picture}(12,5)(-6,-2)
\put(-5,0){\vector(1,0){10}}
\put(0,0){\vector(0,1){0.5}}
\put(0,0){\vector(0,-1){0.5}}
\put(3.5,0){\vector(0,1){0.5}}
\put(3.5,0){\vector(0,-1){0.5}}
\put(-3.5,0){\vector(0,1){0.5}}
\put(-3.5,0){\vector(0,-1){0.5}}
\put(0,0){\ellipse{9.5}{4.5}}
\put(-4.7,-0.1){\vector(0,1){0.1}}
\put(-2,0){\circle*{0.1}}
\put(2,0){\circle*{0.1}}
\put(-2,0){\circle{2}}
\put(2,0){\circle{2}}
\put(-2.95,0.1){\vector(0,-1){0.1}}
\put(1.05,0.1){\vector(0,-1){0.1}}
\put(-3.5,0.8){\mbox{$\tilde\mu_1$}}
\put(0,0.8){\mbox{$\tilde\mu_2$}}
\put(3.5,0.8){\mbox{$\tilde\mu_3$}}
\put(-2,1.2){\mbox{$T_1$}}
\put(2,1.2){\mbox{$T_2$}}
\put(-3,2.2){\mbox{$T_3$}}
\put(-2.1,-0.4){\mbox{$0$}}
\put(+1.9,-0.4){\mbox{$1$}}
\put(4,2.5){\mbox{$\infty$}}
\end{picture}
\end{center}
\caption{Reflections $\tilde\mu_1$, $\tilde\mu_2$ and $\tilde\mu_3$.}
\label{fig:refl}
\end{figure}

As \eqref{eq:pos-sing}, we identify $S^2_{p_1,p_2,p_3}$
with $\C\setminus\{0,1\}$.
Let $\mu$ be the reflection (i.e., 
conformal transformations reversing orientation) on 
$\C\setminus\{0,1\}$
along the real axis.
Then there exist reflections $\tilde \mu_k$ 
on the universal cover $\tilde S^2_{p_1,p_2,p_3}$ such that 
\[
   \pi\circ\tilde\mu_k=\mu\circ\pi,\qquad (k=1,2,3),
\]
and
the deck transformations  $T_j$ ($j=1, 2, 3$) are represented as
 (see Figure~\ref{fig:refl})
\begin{equation}\label{eq:deck}
  T_1 = \tilde\mu_1\circ\tilde\mu_2, \qquad
  T_2 = \tilde\mu_2\circ\tilde\mu_3, \qquad
  \mbox{and}
  \qquad
  T_3 = \tilde \mu_3\circ \tilde \mu_1.
\end{equation}
Let $\tilde F$ be a solution of the equation \eqref{eq:ode}.
Then $\overline{\tilde F\circ\tilde\mu_k}$ ($k=1, 2, 3$)
is also a solution of \eqref{eq:ode} because of 
\begin{equation}\label{eq:GQ-symm}
   \overline{Q\circ \mu} =  Q,\qquad\mbox{and}
   \qquad
   \overline{G\circ \mu} =  G.
\end{equation}
Hence, 
there exist matrices $\rho_{\tilde F}(\tilde\mu_k)\in \SL(2,\C)$ such that
\[
   \overline{\tilde F\circ \tilde\mu_k}=
     \tilde F\cdot \rho_{\tilde F}(\tilde\mu_k) \qquad (k=1,2,3).
\]
Since $\tilde\mu_k\circ\tilde\mu_k=\id$, it hold that
\begin{equation}\label{eq:square}
   \rho_{\tilde F}(\tilde\mu_k)\cdot
   \overline{\rho_{\tilde F}(\tilde\mu_k)}=\id\qquad
     (k=1,2,3),
\end{equation}
and by \eqref{eq:deck}, we have
\begin{equation}\label{eq:deck-mu}
  \begin{array}{r}
   \rho_{\tilde F}(\tilde\mu_1)\cdot\overline{\rho_{\tilde F}(\tilde\mu_2)} 
       = \rho_{\tilde F}(T_1),\quad
   \rho_{\tilde F}(\tilde\mu_2)\cdot\overline{\rho_{\tilde F}(\tilde\mu_3)} 
       = \rho_{\tilde F}(T_2),\hspace{15mm}\\[8pt]
   \quad\mbox{and}\quad
   \rho_{\tilde F}(\tilde\mu_3)\cdot\overline{\rho_{\tilde F}(\tilde\mu_1)} 
     = \rho_{\tilde F}(T_3).
  \end{array}
\end{equation}

If there exists a solution $\tilde F$ of \eqref{eq:ode} such that 
\begin{equation}\label{eq:mu-unitary}
   \rho_{\tilde F}(\tilde\mu_k)\in\SU(2) \qquad (k=1,2,3),
\end{equation}
it follows that $\rho_{\tilde F}(T_j)\in\SU(2)$ for $j=1, 2, 3$
because of \eqref{eq:deck-mu}.
For such an $\tilde F$, we set $g$ as in \eqref{eq:g-def} and 
$d\sigma^2$ as in \eqref{eq:pullback}.
Then $\rho_g(T)\in\SU(2)$ for each deck transformation $T$ on
$\tilde S_{p_1,p_2,p_3}$.
This implies that $d\sigma^2\in\Met(S^2)$.
Moreover, the divisor of $d\sigma^2$ is $D$
by \eqref{S-rel}.

Thus, it is enough to show 
that there exists a solution $\tilde F$ of \eqref{eq:ode} which satisfy
\eqref{eq:mu-unitary}.
To do this, we use the following argument similar to the proof of 
Proposition~6.7 and Proposition~5.6 in \cite{ruy1}.

\step{}    \label{step:1}
  Let $z_0$ be a point on the segment $(-\infty,0)$ of the real axis on 
  $\C$ and take a point $\tilde z_0\in\tilde S_{p_1,p_2,p_3}$ such
  that $\pi(\tilde z_0)=z_0$
  and $\tilde \mu_1(\tilde z_0)=\tilde z_0$.
  Take the solution $\tilde F$ of \eqref{eq:ode} with
  an initial condition $\tilde F(\tilde z_0)=\id$.
  By \eqref{eq:GQ-symm}, $\overline {\tilde F\circ\tilde\mu_1}$ 
  is also a solution of \eqref{eq:ode}.
  Moreover, since $\tilde \mu_1(\tilde z_0)=\tilde z_0$,
  $\overline{\tilde F\circ\tilde\mu_1}$ has the same initial 
  condition as $\tilde F$.
  Hence we have $\overline{\tilde F\circ\tilde\mu_1}=\tilde F$, and then
  $\rho_{\tilde F}(\tilde\mu_1)=\id$.
\step{} \label{step:2}
  Let $\tilde F$ be as in the previous step.
  By \eqref{eq:deck-mu} and Lemma~\ref{lem:eigen}, 
  the eigenvalues of $\rho_{\tilde F}(\tilde\mu_2)$ are $-e^{\pm iB_2}$.
  Moreover, $\sin B_2\neq 0$ because $\beta_2$ is not an integer.
  In particular, $\rho_{\tilde F}(\tilde\mu_2)$ is semisimple.
  By Fact \ref{fact3} and \eqref{eq:square},
  there exists a matrix $u\in \SL(2,\R)$ such that
  \[
       u^{-1}\cdot \rho_{\tilde F}(\tilde\mu_2)\cdot u =
       \pmatrix{
         -e^{iB_2} & 0 \cr
          0  & -e^{-iB_2}}.
  \]
  Let $\hat F := \tilde F\cdot  u$.
  Then $\hat F$ is also a solution of \eqref{eq:ode} and
  \begin{equation}\label{eq:f-hat-rho}
   \begin{array}{rl}
      \rho_{\hat F}(\tilde\mu_1) &= 
               u^{-1}\cdot \rho_{\tilde F}(\tilde\mu_1)\cdot  \bar {u }=
               u^{-1}\cdot \rho_{\tilde F}(\tilde\mu_1)\cdot  u =\id,\\
      \rho_{\hat F}(\tilde\mu_2) &= 
               u^{-1}\cdot \rho_{\tilde F}(\tilde\mu_2)\cdot  \bar {u }=
               u^{-1}\cdot \rho_{\tilde F}(\tilde\mu_2)\cdot  u =
             \displaystyle\pmatrix{
               -e^{iB_2} & 0 \cr
               0  & -e^{-iB_2}}
    \end{array}
  \end{equation}
  hold because $u$ is a real matrix.
\step{} \label{step:3}
  Let $\hat F$ be as in Step~\ref{step:2}.
  By Fact \ref{fact3} and \eqref{eq:square},
  $\rho_{\hat F}(\tilde\mu_3)$ can be written as
  \[
      \rho_{\hat F}(\tilde\mu_3)=\pmatrix{
                        q & i \delta_1 \cr
                        i\delta_2 & \bar q},
      \qquad (\delta_1,\delta_2\in \R, \quad q\bar q+\delta_1\delta_2=1).
  \]  
  Then by \eqref{eq:deck-mu} and Lemma~\ref{lem:eigen}, we have
  \[
     q =  \frac{i}{\sin B_1}(\cos B_2 + e^{iB_1}\cos B_3).
  \]
  Hence, by the assumption \eqref{eq:condition},
  \begin{eqnarray*}
   \delta_1\delta_2 &=& 1-q \bar q\\
      &=&1-
        \frac{\cos^2 B_1+\cos^2 B_2+\cos^2 B_3+2\cos B_1\cos B_2\cos B_3}{
               \sin^2 B_1}>0.
  \end{eqnarray*}
  Let 
  \[
      \check F = \hat F\cdot \pmatrix{
                     (\delta_1/\delta_2)^{1/4}& 0 \cr
                        0       & (\delta_2/\delta_1)^{1/4}}
                    \in\SL(2,\R).
  \]
  Then $\check F$ is a solution of \eqref{eq:ode} and
  \[
    \rho_{\check F}(\tilde\mu_1)=\id,\quad
    \rho_{\check F}(\tilde\mu_2)=\pmatrix{
                     -e^{iB_1} & 0 \cr
                     0  & -e^{-iB_1}},
   \quad
   \rho_{\check F}(\tilde\mu_3)=\pmatrix{
                     q & i\delta \cr
                     i\delta  & \bar q}
   \in\SU(2),
  \]
  where $\delta=(\delta_1\delta_2)^{1/2}$.
  Thus, we have a desired metric $d\sigma^2$ induced from
  $g=-d\check F_{12}/d\check F_{11}$.
\end{Proof*}

Metrics with conical singularities are closely related to 
\cmcone{} (constant mean curvature $1$) surfaces in the hyperbolic
$3$-space $H^3$.
In fact, as shown in \cite[Theorem~2.2]{UY3}
the set $\Met(\Sigma)$ corresponds
bijectively to the set of branched \cmcone{} immersions of $\Sigma$
excluded finite points of finite total curvature 
with prescribed hyperbolic Gauss map.
The one direction of the correspondence is given as follows:
Let $x:M:=\Sigma\setminus\{p_1,\ldots,p_n\}\to H^3$ be a conformal 
\cmcone{} immersion whose induced metric $ds^2$ is complete and of 
finite total curvature.
We set
\[
  d\sigma^2_x:=(-K)ds^2,
\]
where $K$ is the Gaussian curvature of $ds^2$.
Then it can be extended to a pseudometric on $\Sigma$ and
$d\sigma^2_x\in \Met(\Sigma)$ holds (cf.~\cite{B}). 
The converse correspondence is described in \cite[Section~2]{UY3}.

\begin{Def}
  A regular end of a \cmcone{}
   immersion $x$ is called Type I,
  if the Hopf differential of the immersion has pole of order $2$
  at the end (cf.\ \cite{ruy2}).
\end{Def}

As seen in \cite[Section 5]{UY1}, 
{\it a regular end of the \cmcone{} immersion $x$ is
   asymptotic to a certain catenoid cousin end if and only if
   the end is of Type I and embedded.}

Using the same argument of the proof of 
Theorem~\ref{thm:irred}, we can classify 
the set of irreducible \cmcone{} surface 
in the hyperbolic $3$-space of genus zero,
with three ends asymptotic the catenoid cousins.

\begin{Thm}\label{prop:surface}
  If a triple of real numbers $(\beta_1,\beta_2,\beta_3)$ satisfy
  \eqref{eq:condition} for $B_j=\pi(\beta_j+1)$ $(j=1, 2, 3)$,
  there exists a unique irreducible conformal \cmcone{} 
  immersion $x:S^2_{p_1,p_2,p_3}\to H^3$
  such that all ends $p_1,p_2,p_3$ are of embedded and Type I,
  and the order of the pseudometric $d\sigma^2_x$ at $p_j$ is $\beta_j$.
  Conversely, 
  any conformal irreducible immersed \cmcone{} surface of 
  genus zero with three embedded Type I ends are obtained in such a manner.
\end{Thm}

\begin{Proof}
  In fact, such a surface is realized by an conformal \cmcone{} immersion
  $x:\C\setminus\{0,1\}\to H^3$ with the following properties:
  \begin{enumerate}
     \item\label{item:surf:1} 
           Since three ends are of Type I, the Hopf 
           differential $Q$ of $x$ has poles of order $2$ 
           at the ends $0$, $1$ and $\infty$.
           Then necessarily it has two zeros $q_1$, $q_2$ of order $1$ 
           on $\C\setminus\{0,1\}$.
            (When $Q$ has the only zero 
                of order two, we regard
            $q_1=q_2$.)
     \item\label{item:surf:2} 
           Since all ends are regular and embedded, the hyperbolic 
           Gauss map $G$ of $f$ has two branch point 
           of order $1$ at zeroes of $Q$, and
           no branch point elsewhere (see \cite{UY4}).
           (Hence, $G$ is a meromorphic function of degree $2$.) 
     \item\label{item:surf:3} 
           $\tilde S(d\sigma^2_x)-S(G)=2Q$ holds (see \cite[(2.3)]{UY3}).
     \item\label{item:surf:4} 
           Since the order of the pseudometric $d\sigma^2_x$
           at the umbilic points $q_1,q_2$ are equal to the order of 
           zeroes
            of $Q$ (see \cite{UY4}), 
           $d\sigma^2_x$ has the divisor of the form
           \[
              D':=\beta_1 p_1 +\beta_2 p_2 + \beta_3 p_3+ q_1 +  q_2,\qquad
                     (\beta_j>-1).
           \]
           The metric $d\sigma^2_x$ is irreducible if and only if 
           \eqref{eq:condition} holds (cf.~Appendix~\ref{app:A}).
  \end{enumerate}
  We may assume that, $p_1=0$, $p_2=1$ and $p_3=\infty$.
  By \ref{item:surf:1}, \ref{item:surf:2} and \ref{item:surf:3}, 
  the top term of $2Q$ of Laurent expansion at $z=p_j$ is the
  same as that of $\tilde S(d\sigma^2_x)$.
  Thus we have
\[
  Q=\frac{1}{2}
  \left(-\frac{\beta_j(\beta_j+2)}2\frac1{(z-p_j)^2}+\cdots\right)dz^2\qquad
  (j=1,2,3).
\]
Since the Hopf differential $Q$ is holomorphic on $S^2_{p_1,p_2,p_3}$,
we have
\begin{equation}\label{eq:Q-surface}
  Q=\frac{1}{2}\left(
               \frac{c_3 z^2+(c_2-c_1-c_3)z+c_1}{z^2(z-1)^2}
            \right)\,dz^2,
\end{equation}
where $c_j=-\beta_j(\beta_j+2)/2\in\R$  ($j=1$, $2$, $3$).
By (2), $G^*d\sigma^2_0$ has the divisor of the form
\[
   D_G:=q_1 +  q_2.
\]
Since the hyperbolic Gauss map $G$ has
an ambiguity of M\"obius transformations,
we can set
\begin{equation}\label{eq:G-surface}
      G =  z + \displaystyle\frac{(q_1-q_2)^2}{2\{2z-(q_1+q_2)\}}.
\end{equation}
Then by \cite[Theorem~3.1]{ruy1},
we can see the uniqueness of  
irreducible \cmcone{} immersion $x$ with the hyperbolic
Gauss map $G$ and the Hopf differential $Q$.
So it is sufficient to show the existence of such a surface.
The following proof is almost same as that of Theorem~\ref{thm:irred}:
Let $\mu$ be the reflection 
with respect to the real axis and take
the reflections $\tilde \mu_j$ ($j=1,2,3$) 
on the universal cover $\tilde S^2_{p_1,p_2,p_3}$
as in the proof
of Theorem~\ref{thm:irred} (see Figure~\ref{fig:refl}).

Let $\tilde F$ be a solution of the equation 
\begin{equation}\label{eq:ODE2}
  d\tilde F \cdot \tilde F^{-1}
        =\pmatrix{G& -G^2 \cr 1 & -G\hphantom{^2}}\frac{Q}{dG}
\end{equation}
for $G$ in \eqref{eq:G-surface} and $Q$ in \eqref{eq:Q-surface}.
Then $\overline{\tilde F\circ\tilde\mu_k}$ ($k=1, 2, 3$)
is also a solution of \eqref{eq:ode} because of 
\begin{equation}\label{eq:GQ-symm2}
   \overline{Q\circ \mu} =  Q,\qquad
   \mbox{and}%
   \qquad
   \overline{G\circ \mu} =  G.
\end{equation}
Hence, there exist matrices 
 $\rho_{\tilde F}(\tilde\mu_j)\in \SL(2,\C)$
 such that
\begin{equation}
   \overline{\tilde F\circ \tilde\mu_j}
   =\tilde F\cdot \rho_{\tilde F}(\tilde\mu_j) 
\qquad (j=1,2,3).
\end{equation}
Now by the completely same argument as in the proof of Theorem \ref{thm:irred},
there exists a solution $\tilde F$ of \eqref{eq:ODE2} such that 
$\rho_{\tilde F}(T_j)\in\SU(2)$ for $j=1, 2, 3$.
For such an $\tilde F$, we set $g$ as in \eqref{eq:g-def} and 
$d\sigma^2$ as in \eqref{eq:pullback}.
Then $d\sigma^2\in\Met(S^2)$ has the divisor $D'$.
By \cite[Theorem~2.2]{UY3},
there exists a branched \cmcone{} immersion 
$x:S^2_{p_1,p_2,p_3}\to H^3$
whose  hyperbolic Gauss map and Hopf differential are $G$ and $Q$
respectively such that $d\sigma^2_x=d\sigma^2$.
One can easily check that the metric given by
\[
   ds^2{}^\sharp:=(1+|G|^2)^2 \left|\frac{Q}{dG}\right|^2
\]
is positive definite and complete.
Thus by \cite[Lemma~2.3]{ruy1}, so is the 
first fundamental form $ds^2$.
Hence $x$ is the desired one. 
\end{Proof}

\begin{Prop}\label{prop:ta}
  Let $x:S^2_{p_1,p_2,p_3}\to H^3$ be complete \cmcone{} surface
  with three ends of Type I.
  Then the total absolute curvature $\TA$ of $x$
  is greater than or equal to $4\pi$.
\end{Prop}

\begin{Rmk}
In \cite{UY1}, the authors showed 
$\TA>2\pi$ for three ended \cmcone{} surfaces.
The estimate in Proposition~\ref{prop:ta} is sharper than this.
\end{Rmk}

\begin{Proof}
 The associated pseudometric $d\sigma^2_x$ has
 the divisor $D':=\beta_1p_1+\beta_2p_2+\beta_3p_3+q_1+q_2$
 where $q_1,q_2$ are umbilic points of $x$.
 Then we have
 \begin{eqnarray*}
    \frac{1}{2\pi}\TA&=&\frac1{2\pi}\int_{S^2_{p_1,p_2,p_3}}(-K)\,dA_{ds^2}=
        \frac1{2\pi}\int_{S^2_{p_1,p_2,p_3}}\,dA_{d\sigma^2}\\
         &=&\chi(S^2)+|D|=4+\beta_1+\beta_2+\beta_3.
 \end{eqnarray*}
 On the other hand, the metric $d\sigma^2_x$ induces
 a monodromy representation 
 $\rho_g:\pi_1(S^2_{p_1,p_2,p_3})\to \PSU(2)$.
 As seen in the proof of Theorem \ref{thm:irred}, 
 it can be lifted to a representation
 $\rho_{\tilde F_g}:\pi_1(S^2_{p_1,p_2,p_3})\to \SU(2)$.
Let $T_j$  ($j=1,2,3$) be the deck transformation 
corresponding the the loop surrounding $p_j$.
Then the eigenvalues of $\rho_{\tilde F_g}$
is $-e^{\pm iB_j}$ where $B_j:=\pi(\beta_j+1)$, which can be
proved by the same method as in Lemma \ref{lem:eigen}.
By Lemma~\ref{lem:A} in Appendix~\ref{app:A},
we have
\[
  \cos^2 B_1 +\cos^2 B_2 + \cos^2 B_3+
               2\cos B_1\cos B_2 \cos B_3<1.
\]
Then it can be easily seen that
$B_1+B_2+B_3\ge \pi$,
which yields $\beta_1+\beta_2+\beta_3\ge -2$.
\end{Proof}
\section{Reducible metrics with three singularities.}
\label{sec:reducible}

In this section, we give a
necessary and sufficient 
condition for the existence of reducible metrics
with given divisors.
As in the previous section, we identify $S^2=\C\cup\{\infty\}$, and 
set $(p_1,p_2,p_3)=(0,1,\infty)$.

\begin{Lemma}\label{lem:integer}
  Let $d\sigma^2\in\Met(S^2)$ be an 
  $\Hyp^3$-reducible $($resp.\ $\Hyp^1$-reducible$)$ pseudometric
  with divisor $D$ as in \eqref{eq:D}.
  Then all of $\beta_j$'s are integers $($resp.\ 
  exactly one of $\beta_j$'s is an integer$)$.
  Namely, at least one of $\beta_j$ is an integer
  for reducible metrics.
\end{Lemma}
\begin{Proof}
  Let $d\sigma^2\in\Met(S^2)$ be an $\Hyp^3$-reducible metric
  with divisor $D$ as in \eqref{eq:D}.
  Then, by definition, the representation $\rho_g$ as in 
  \eqref{eq:g-repr} is trivial.
  Hence $g$ satisfying \eqref{eq:pullback} is a single-valued
  meromorphic function on $S^2_{p_1,p_2,p_3}$.
  Since $p_j$ ($j=1, 2, 3$) are conical singularities of
  $d\sigma^2$, $g$ can be extended to a meromorphic function on 
  $S^2$.
  Hence $\beta_j$ ($j=1, 2, 3$) are integers.

\medskip

  Next, assume $d\sigma^2$ is $\Hyp^1$-reducible and all $\beta_j$'s are 
  non-integral numbers.
  Since $d\sigma^2$ is reducible, we can choose $g$ as in
  \eqref{eq:pullback} such that $\rho_g(T)$ are diagonal for
  all $T\in\pi_1(S^2_{p_1,p_2,p_3})$.
  Then $g\circ T_k=e^{2\pi i(\beta_k+1)} g$ because of
  \eqref{eq:g-expand}.
  Hence $g_1:=z^{-\beta_1-1} g$
  is single-valued on 
  $\tilde S^2_{p_2,p_3}$, and 
  $g_2:=(z-1)^{-\beta_2-1} g_1$
  is
  single-valued on $\tilde S^2_{p_3}$.
  Since $S^2_{p_3}$ is simply-connected, $g_2$
  is single-valued on $S^2$
  and $g$ can be written as
  \begin{equation}\label{eq:g-3-nonint}
    g = z^{\mu}(z-1)^{\nu} \frac{a(z)}{b(z)}\qquad (\mu,\nu\in\R\setminus\Z),
  \end{equation}
  where $a(z)$ and $b(z)$ are mutually prime polynomials whose roots 
  are distinct from $0$ and $1$.
  Then we have
  \begin{equation}\label{eq:dg-3-nonint}
  \begin{array}{rl}
     &dg= z^{\mu-1}(z-1)^{\nu-1} \displaystyle\frac{p(z)}{q(z)}\,dz,\\
     &\quad p(z):=
     \{\nu z + \mu (z-1) \}a(z)b(z)+z(z-1)\{a'(z)b(z)-a(z)b'(z)\},\\
     &\quad q(z):=\{b(z)\}^2.
  \end{array}
  \end{equation}
  Since $p(0)=-\mu \, a(0)b(0)$, $p(1)=\nu\,  a(0)b(0)$, $q(0)=\{b(0)\}^2$ 
  and $q(1)=\{b(1)\}^2$ are not equal to $0$,
  the roots of $p$ and $q$ are distinct from $0$ and $1$.
  Moreover, $p$ and $q$ are mutually prime.
  In fact, assume there exists a common root $\xi$ of $p$ and $q$.
  Then $b(\xi)=0$, and by assumption, $\xi\neq 0,1$ and $a(\xi)\neq 0$.
  Then $0=p(\xi)=\xi(\xi-1)a(\xi)b'(\xi)$ implies $b'(\xi)=0$.
  Hence $\xi$ is a multiple root of $b$.
  Let $b(z)=(z-\xi)^m \tilde b(z)$, where 
  $m\geq 2$ be an integer and $\tilde b(z)$ is a polynomial such that 
  $\tilde b(\xi)\neq 0$.
  Then we have
  \[
      dg = z^{\mu-1}(z-1)^{\nu-1}
         \frac{(z-\xi)r(z) - z(z-1)\tilde b(z)}{
                    (z-\xi)^{m+1}\tilde b(z)^2}\,dz,
  \]
  where $r(z)=(\nu z+\mu(z+1))a\tilde b+z(z-1)(a'\tilde b-a\tilde b')$
  is a polynomial in $z$.
  Since $m\geq 2$, $\xi$ is a ramification point of $g$, and then,
  $d\sigma^2$ has a conical singularity at $\xi$, 
  which is a contradiction.

  Thus $p$ and $q$ are mutually prime whose roots are distinct from $0$
  and $1$.
  If $p$ has a root $\eta$, then $\eta$ is a ramification point of $g$,
  and then, a conical singularity of $d\sigma^2$.
  Hence $p(z)$ must be a constant.
  By \eqref{eq:dg-3-nonint}, $p$ is formally a polynomial of degree 
  $\deg a+\deg b +1$.
  Then the highest term must vanish:
  $\mu+\nu+\deg a - \deg b = 0$.
  This shows that the order of $dg$ at $z=\infty$ must be integer,
  and then $\beta_3$ is an integer.
  This is contradiction, 
  and hence at least one of $\beta_j$'s must be
  integer.

\medskip

  On the other hand, 
  let $d\sigma^2\in\Met(S^2)$ be  $\Hyp^1$-reducible 
  and
  exactly one of  $\beta_j$'s is a non-integral number.
  Without loss of generality, we assume $\beta_1$ is non-integral.
  Take $g$ as in \eqref{eq:pullback}.
  Since $\beta_2$ and $\beta_3$ are integers,
  $g$ is well-defined on 
  the universal cover of $S^2_{p_1}$.
  Here $S^2_{p_1}$ is simply connected.
  Then $g$ is single-valued on $S^2$ itself, and hence, $g$ is a 
  meromorphic function on $S^2$.
  This shows that $\beta_1$ is an integer, contradiction.

\medskip
  
  Hence, if $d\sigma^2\in\Met(S^2)$ with divisor $D$ in \eqref{eq:D}
  is $\Hyp^1$-reducible, exactly one of $\beta_j$'s is an integer.
\end{Proof}

As an immediate consequence, we have the following corollary.

\begin{Cor}\label{added}
Let $d\sigma^2\in \Met(S^2)$ has exactly three singularities
with orders $\beta_1$, $\beta_2$ and $\beta_3$.
Then the following three assertions are true.
\begin{enumerate}
\item $d\sigma^2$ is $\Hyp^3$-reducible if and only if 
all of $\beta_j$'s are integers.
\item 
$d\sigma^2$ is $\Hyp^1$-reducible if and only if
 exactly one of $\beta_j$'s is an integer.
\item
$d\sigma^2$ is irreducible if and only if
all of $\beta_j$'s are non-integers.
\end{enumerate}
\end{Cor}

\begin{Rmk}
For a metric in $\Met(S^2)$ with more than three singularities,
such a simple criterion for reducibility is not expected:
There exists a reducible metric $d\sigma^2\in\Met(S^2)$ with
  divisor
  \[
     D' =\beta_1 p_1 +\beta_2 p_2 + \beta_3 p_3+ q_1 +  q_2
  \]
  such that all $\beta_j$'s are non-integral numbers.
  In fact, 
  \[
     g=c\, z^{\mu}\, (z-1)^{\nu} (z-a)
      \qquad (c\in \C\setminus\{0\},\ 
              a\in \C\setminus\{0,1\})
  \]
  induces such a metric whenever $\mu+\nu$ is not an integer.
On the other hand, we can construct an
irreducible metric with divisor $D'$ such that
$\beta_1,\beta_2,\beta_3\not\in \mathbf Z$:
The metric $d\sigma_x^2\in\Met(S^2)$ obtained in
Theorem~\ref{prop:surface} is the desired one.
\end{Rmk}

\subsection*{$\Hyp^3$-reducible case:}
  First, we consider the case of $\Hyp^3$-reducible.
  In this case, $\beta_1$, $\beta_2$ and $\beta_3$ are integers and 
  $g$ in  \eqref{eq:pullback} is single-valued on $S^2$, 
  i.e., a rational function on $\C\cup\{\infty\}$.

  Without loss of generality, we  assume
\begin{equation}\label{eq:order}
   \beta_1 \leq \beta_2\leq \beta_3.
\end{equation}
  Let $g$ be a rational function such that $d\sigma^2$ as in
  \eqref{eq:pullback} with the divisor $D$ in \eqref{eq:D} satisfying
  \eqref{eq:order}.
  Then the ramification points of $g$ are $0$, $1$ and $\infty$ whose 
  orders are $\beta_1$, $\beta_2$ and $\beta_3$ respectively.
  By the Riemann-Hurwicz formula, 
\[
    \deg g=\frac{1}{2}(\beta_1+\beta_2+\beta_3)+1 
          \leq \frac{1}{2}\beta_1 + \beta_3 + 1
          < (\beta_1+1)+(\beta_3+1)
\]
  holds.
  Then we have $g(p_1)\neq g(p_3)$, and similarly, $g(p_2)\neq g(p_3)$.
  Thus, by a suitable change as \eqref{eq:g-change}, we may assume
  $g(p_1)=g(0)\neq \infty$, $g(p_2)=g(1)\neq \infty$, 
  and $g(p_3)=g(\infty)=\infty$.
  Under these assumptions, we can write
\begin{equation}\label{eq:dg-int}
    dg = c
     \frac{z^{\beta_1}(z-1)^{\beta_2}}{\prod_{j=1}^N(z-a_j)^2}\,dz,\qquad
   \beta_3 = \beta_1+\beta_2-2N,
\end{equation}
  where $c\neq 0$ is a constant, and $a_1$, \ldots,
  $a_N\in\C\setminus\{0,1\}$ are mutually distinct numbers.
  Conversely, if there exists $g$ which satisfies
  \eqref{eq:dg-int}, then
  we have $d\sigma^2\in\Met(S^2)$ with desired singularities.
  Computing residues at $z=a_1$, \ldots, $a_N$, we have 
\begin{Prop}\label{prop:reducible1}
  Let $\beta_1$, $\beta_2$ and $\beta_3$ are positive integers
  satisfying \eqref{eq:order}.
  Then there exists $d\sigma^2\in\Met(S^2)$ with divisor $D$ in \eqref{eq:D}
  if and only if 
  there exists a non-negative integer $N$ and mutually distinct complex
  numbers $a_1$, \ldots, $a_N\in\C\setminus\{0,1\}$ which satisfy
  \begin{equation}\label{eq:dg-int-N}
    \beta_{1}+\beta_{2}-\beta_{3} = 2N 
  \end{equation}
  and
  \begin{equation}\label{eq:dg-int-res}
    \frac{\beta_{1}}{a_j}+
              \frac{\beta_{2}}{a_j-1}-
              \sum_{k\neq j}
              \frac{2}{a_j-a_k} =0,\qquad
             (j=1,\ldots, N).
  \end{equation}
\end{Prop}
\begin{Rmk}
  An $\Hyp^3$-reducible
  metric admits a three parameter space of 
  deformation as in \eqref{eq:deform-sp} which preserves 
  the divisor 
  $D$ and the Schwarzian derivative.
  For a given triple $(\beta_1,\beta_2,\beta_3)$  as in 
  Proposition~\ref{prop:reducible1}, such a deformation space is
  determined uniquely.
  In fact, assume two metrics $d\sigma_1^2=g_1^*d\sigma_0^2$ and 
  $d\sigma_2^2=g_2^*d\sigma_0^2$ have the same divisor $D$.
  Then both $\tilde S(d\sigma_1^2)=S(g_1)$ and $\tilde S(d\sigma_2^2)=S(g_2)$ 
  are equal to \eqref{eq:schwarz-2}.
  Then there exists $a\in\PSL(2,\C)$ such that $g_2=a\star g_1$.
\end{Rmk}

  Observing the equation \eqref{eq:dg-int-res} for small $N$,
  one can easily see the following facts:
  If $N$ in \eqref{eq:dg-int-N} is $0$, trivially the metric exists.

  When $N=1$, \eqref{eq:dg-int-res} has the unique solution
  $a_1=\beta_1/(\beta_1+\beta_2)\neq 0$, $1$.

  Assume $N=2$.
  Then $\beta_3=\beta_1+\beta_2-4$.
  By \eqref{eq:order}, this implies $4\leq\beta_1\leq\beta_2$.
  In this case, it is easy to show that
  the system of equations \eqref{eq:dg-int-res} has
  the unique solution up to the order of $a_j$'s.

  Hence, we have 
\begin{Cor}
  Let $\beta_1$, $\beta_2$ and $\beta_3$ are positive integers
  satisfying \eqref{eq:order}
  and
  \[
     \beta_1+\beta_2-\beta_3=2N\qquad (N=0,1\mbox{ or }2).
  \]
  Then, there exists an $\Hyp^3$-reducible metric
  $d\sigma^2\in\Met(S^2)$ with the divisor $D$
  as in \eqref{eq:D}.
  Moreover, such a metric is unique up to three parameter deformation as 
  in \eqref{eq:deform-sp}.
\end{Cor}

\subsection*{$\Hyp^1$-reducible case:}
  For $\Hyp^1$-reducible case, one of the $\beta_j$'s must be
  an integer
  because of Lemma~\ref{lem:integer}.
  We  assume that $\beta_1$ and $\beta_3$ are non-integral numbers,
  and $\beta_2$ is a positive integer:
\[
   \beta_1,\,\,\beta_3 \not\in \Z, \qquad 
   \beta_2 \in \Z.
\]
  Then  we can choose $g$ as in \eqref{eq:pullback} such that
  \begin{equation}\label{eq:g-2-nonint}
     g = z^{\mu} \varphi(z),\qquad \mu\in\R\setminus\Z,
  \end{equation}
  where $\varphi(z)$ is a rational function.
  Here, such a normalization of $g$ is unique up to the change
  $g\mapsto t g$ or $g\mapsto t/g$ for 
  each non-zero constant $t$.
  Moreover, 
\begin{equation}\label{eq:red1-change}
  (tg)^* d\sigma_0^2\qquad (t\in\R^+)
\end{equation}
  gives a non-trivial deformation of the metric preserving the divisor.
  This is the one parameter deformation as
  \eqref{eq:deform-sp}.

  Let $g$ be a function as in \eqref{eq:g-2-nonint} such that 
  $d\sigma^2$ as in \eqref{eq:pullback} has the divisor $D$ 
  in \eqref{eq:D}.
  Replacing $g$ with $1/g$, we may assume $\varphi(1)\neq \infty$ without 
  loss of generality.
  Under this assumption, $dg$ has zero of order $\beta_2$ at $z=1$:
\begin{equation}\label{eq:dg-2-nonint}
    dg = c z^{\nu_1}\frac{(z-1)^{\beta_2}}{\prod_{j=1}^{N}(z-a_j)^2}\,dz
   \qquad (\nu_1=\mu-1),
\end{equation}
  where $c\neq 0$ is a constant, and $a_1$,
  \ldots,$a_N\in\C\setminus\{0,1\}$
  are mutually different numbers.

  We denote the order of $dg$ at $z=\infty$ by $\nu_3$:
\begin{equation}\label{eq:nu}
  \nu_3 = -\nu_1-\beta_2+2N -2.
\end{equation}
So the following four cases occur.
\begin{enumerate}
\renewcommand{\labelenumi}{(\alph{enumi})}
\renewcommand{\theenumi}{(\alph{enumi})}
\item\label{item:red1}
      $\nu_1>-1$ and $\nu_3>-1$.
      In this case, $\beta_1=\nu_1$, $\beta_3=\nu_3$.
      Hence $\beta_1+\beta_2+\beta_3=2N-2$.
\item\label{item:red2} 
      $\nu_1<-1$ and $\nu_3<-1$.
      In this case, $\beta_1=-\nu_1-2$, $\beta_3=-\nu_3-2$.
      Hence $\beta_1-\beta_2+\beta_3=-2N-2$.
\item\label{item:red3} 
      $\nu_1>-1$ and $\nu_3<-1$.
      In this case, $\beta_1=\nu_1$, $\beta_3=-\nu_3-2$.
      Hence $\beta_1+\beta_2-\beta_3=2N$.
\item\label{item:red4} 
      $\nu_1<-1$ and $\nu_3>-1$.
      In this case, $\beta_1=-\nu_1-2$, $\beta_3=\nu_3$.
      Hence $\beta_1-\beta_2-\beta_3=-2N$.
\end{enumerate}
For the cases \ref{item:red1} and \ref{item:red3}, 
there exists a meromorphic function $g$ on the universal cover of 
$\C\setminus\{0\}$ satisfying \eqref{eq:dg-2-nonint} if and only if
\begin{equation}\label{eq:red1-eq-1}
      \frac{\beta_1}{a_j}+
               \frac{\beta_2}{a_j-1}-
               \sum_{k\neq j}
               \frac{2}{a_j-a_k} =0,\qquad
                 (j=1,\ldots, N)
\end{equation}
holds, and for the cases \ref{item:red2} and \ref{item:red4}, 
there exists  $g$ satisfying \eqref{eq:dg-2-nonint} if and only if
\begin{equation}\label{eq:red1-eq-2}
      \frac{-\beta_1-2}{a_j}+
               \frac{\beta_2}{a_j-1}-
               \sum_{k\neq j}
               \frac{2}{a_j-a_k} =0,\qquad
                 (j=1,\ldots, N)
\end{equation}
holds.

Then we have
\begin{Prop}\label{prop:reducible2}
  Let $\beta_1$, $\beta_3$ be non-integral real numbers
  grater than $-1$
  and $\beta_2$ a positive integer.
  Then there exists $d\sigma^2\in\Met(S^2)$ with divisor $D$ in \eqref{eq:D}
  if and only if one of the following occurs{\rm :}
  \begin{enumerate}
  \item\label{item:red-a}
    There exists non-negative integer
    $N$ and mutually distinct complex
    numbers $a_1$, \ldots, $a_N\in\C\setminus\{0,1\}$ which satisfy
    $\beta_1+\beta_2+\beta_3=2N-2$ and 
    \eqref{eq:red1-eq-1} holds.
  \item\label{item:red-b}
    There exists non-negative integer $N$ and mutually distinct complex
    numbers $a_1$, \ldots, $a_N\in\C\setminus\{0,1\}$ which satisfy
    $\beta_1-\beta_2+\beta_3=-2N-2$ and 
    \eqref{eq:red1-eq-2} holds.
  \item\label{item:red-c}
    There exists non-negative integer $N$ and mutually distinct complex
    numbers $a_1$, \ldots, $a_N\in\C\setminus\{0,1\}$ which satisfy
    $\beta_1+\beta_2-\beta_3=2N$ and 
    \eqref{eq:red1-eq-1} holds.
  \item\label{item:red-d}
    There exists non-negative integer $N$ and mutually distinct complex
    numbers $a_1$, \ldots, $a_N\in\C\setminus\{0,1\}$ which satisfy
    $\beta_1-\beta_2-\beta_3=-2N$ and 
    \eqref{eq:red1-eq-2} holds.
  \end{enumerate}
  Moreover, such a metric is unique up to the change in \eqref{eq:red1-change}.
\end{Prop}

As well as $\Hyp^3$-reducible case, we classify the metrics for
$N\leq 2$.
It is easy to show the following lemma:
\begin{Lemma}\label{lem:alg}
 Let $m$ and $N$ be a positive integer, and $\nu$ a non-integral real
 number.
 Consider the following equations 
on $a_1$, \ldots, $a_N${\rm :}
\begin{equation}
  \frac{\nu}{a_j}+\frac{m}{a_j-1}-
  \sum_{%
    \begin{array}{c}\scriptscriptstyle k\neq j\\[-6pt]
                    \scriptscriptstyle 1\leq k\leq N\end{array}%
  }\frac{2}{a_j-a_k}=0,\qquad (j=1,\ldots,N).
\end{equation}
 Then
 \begin{enumerate}
  \item If $N=1$, the equation has the unique solution
	$a_1=\nu/(\nu+m)$, which is different to $0$ and $1$.
  \item If $N=2$, the equation has a solution if and only if 
        $m\neq 1$.
        The both solutions $a_1$, $a_2$ are distinct to $0$ and $1$
        if $m\neq 1$.
 \end{enumerate}
\end{Lemma}
Using this Lemma, we have the following existence and non-existence 
result.
\begin{Cor}\label{cor:nonex}
  Let $\beta_2=1$, and $\beta_1$ and  $\beta_3$ $(>-1)$ be
  non-integral 
  real numbers satisfying $\beta_1+\beta_3=1$.
  Then there exists no metric $d\sigma^2\in\Met(S^2)$ with divisor
  $D$ in \eqref{eq:D}.
\end{Cor}
\begin{Proof}
  This is the case \ref{item:red-a} in
  Proposition~\ref{prop:reducible2} for $N=1$.
  Since $\beta_2=1$, the equation \eqref{eq:red1-eq-1} has 
  no solution because of Lemma~\ref{lem:alg}.
\end{Proof}

\begin{Cor}\label{cor:ex1}
  Let $\beta_1$, $\beta_3$
  be  non-integral numbers $(>-1)$
  and $\beta_2$ an integer satisfying one of the following 
  cases for some non-negative integer $n$:
  \begin{enumerate}
   \item $\beta_1+\beta_3=2n-1$, $\beta_2=2n+2N+1$ $(N=0,1,2)$,
   \item $\beta_1+\beta_3=2n$, $\beta_2=2n+2N+2$ $(N=0,1,2)$,
   \item $\beta_1-\beta_3=2n+1$, $\beta_2=2n+2N+1$ $(N=0,1,2)$,
   \item $\beta_1-\beta_3=2n$, $\beta_2=2n+2N$ $(N=0,1,2)$.
  \end{enumerate}
  Then there exists an $\Hyp^1$-reducible metric
  $d\sigma^2\in\Met(S^2)$ with divisor $D$ in \eqref{eq:D}.
\end{Cor}
\begin{Proof}
  The first two cases follows from \ref{item:red-b} of
  Proposition~\ref{prop:reducible2}, and others from 
  \ref{item:red-d} of  Proposition~\ref{prop:reducible2}.
\end{Proof}

\appendix
\section{}\label{app:A}
In this appendix, we shall prove the following

\begin{appLemma}\label{lem:A}
  Let $a_j\in \SU(2)$ $(j=1,2,3)$ be matrices satisfying 
  $a_1\cdot a_2\cdot a_3=\id$.
  Then the following inequality holds{\rm :}
  \begin{equation}\label{eq:app-cond}
     \cos^2 C_1 +\cos^2 C_2 + \cos^2 C_3+
            2\cos C_1\cos C_2 \cos C_3\leq 1,
  \end{equation}
  where $-e^{\pm iC_j}$ are 
  the eigenvalues of the matrices $a_j$ $(j=1, 2, 3)$.
 Moreover, equality holds if and only if $a_1$, $a_2$ and $a_3$ are
  diagonalizable at the same time. 
\end{appLemma}
\begin{Proof}
  We may assume that $a_1$ is diagonal without loss of generality.
  Then we may set
  \[
    a_1=\pmatrix{
              -e^{iC_1}& 0 \cr
              0 & -e^{- i C_1}
        },
   \qquad 
    a_2=\pmatrix{
                              p & -\bar q \cr
                              q & \hphantom{-}\bar p}
         \quad (|p|^2+|q|^2=1).
  \]
  Since $a_3^{-1}$ has the same eigenvalues as $a_3$,
  we have
  \[
     -2 \cos C_3=\mbox{Tr}(a_1\cdot a_2)=-pe^{iC_1}-\bar pe^{-iC_1}.
  \]
  On the other hand, we have
  \[
      \tr a_2 =       p+\bar p=-2\cos C_2.   
  \]
  Using the above two equations, we have
  \[
    p=\pm i\frac{e^{\pm iC_1} \cos C_2+\cos C_3}{\sin C_1}.
  \]
  Thus we have
  \[
       1\ge |p|^2=\frac1{\sin^2 C_1}
                (\cos^2 C_2 + \cos^2 C_3+
                       2\cos C_1\cos C_2 \cos C_3).
  \]
  This is equivalent to the desired inequality.
  Obviously, equality holds if and only if $|p|=1$, that is, $a_2$
  is diagonal.
\end{Proof}

\section{}\label{app:B}

This is the appendix in \cite{ruy1}.
We attach it here for the sake of convenience. 

Let $\Gamma$ be a subgroup of $\PSU(2)$.

In this appendix, we prove a property of
a set of groups conjugate to $\Gamma$ in $\PSL(2,\C)$ defined by
\[
  C_{\Gamma}:=
  \{  \sigma\in \PSL(2,{\bf C})\,|\,
      \sigma\!\cdot\! \Gamma\! \cdot \!\sigma^{-1}\subset \PSU(2)
  \}.
\]
The authors  wish to thank Hiroyuki Tasaki for 
valuable comments on the first draft of the appendix.

If $\sigma\in C_{\Gamma}$, it is obvious that
$a\!\cdot\! \sigma\in C_{\Gamma}$
for all $a\in \PSU(2)$.
So if we consider the quotient space
\[
   I_{\Gamma}:=C_{\Gamma}/\PSU(2),
\]
the structure of the set $C_{\Gamma}$ is completely determined.
Define a map $\tilde\phi:C_{\Gamma}\to \Hyp^3$ 
by
\[
   \tilde\phi(\sigma):=\sigma^*\!\!\cdot\! \sigma,
\]
where $\Hyp^3$ is the hyperbolic 3-space defined by
$\Hyp^3:=\{a\!\cdot a^*\,|\, a\in \PSL(2,\C)\}$.
Then it induces an
injective map $\phi: I_{\Gamma}\to \Hyp^3$
such that
$\phi\circ \pi=\tilde \phi$, where $\pi:C_{\Gamma}\to I_{\Gamma}$
is the canonical projection.
So we can identify $I_{\Gamma}$ with a subset
$\phi(I_{\Gamma})=\tilde\phi(C_{\Gamma})$
of the hyperbolic 3-space $\Hyp^3$.
The following assertion holds.

\begin{appLemma}\label{lem:B}
  The subset $\phi(I_{\Gamma})$ 
  is a point, a geodesic line,
  or all $\Hyp^3$.
\end{appLemma}

\begin{Proof}
  For each $\gamma\in \Gamma$, we set
  \[
      C_\gamma:=\{\sigma\in \PSL(2,\C)\,|\,
              \sigma\! \cdot\! \gamma \!\cdot\! \sigma^{-1}\in \PSU(2)
        \}.
  \]
  Then we have
  \begin{equation}\label{app:1}
      C_\Gamma:=\bigcap_{\gamma\in \Gamma} C_\gamma.
  \end{equation}
  The condition $\sigma\!\cdot\! \gamma\!\cdot\!  \sigma^{-1}\in \PSU(2)$
  is rewritten as 
  $\sigma^*\!\!\cdot\! \sigma\!\cdot\! \gamma
        =\gamma\!\cdot\! \sigma^*\!\!\cdot\! \sigma$.
  So we have
  \begin{equation}\label{app:2}
     \tilde \phi (C_\gamma)=\Hyp^3 \cap Z_\gamma,
  \end{equation}
  where $Z_\gamma$ is the center of $\gamma\in \Gamma$.

  Assume $\gamma\neq \pm \id$.
  If $\gamma$ is a diagonal matrix, it can easily be checked that
  $Z_{\gamma}$ consists of diagonal matrices in $\PSL(2,\C)$.
  Since any $\gamma\in\Gamma$ can be diagonalized by a matrix in 
  $\SU(2)$, we have $Z_{\gamma}=\{\exp( zT)\,|\,z\in\C\}$, where
  $T\in\Su(2)$ is so chosen that $\gamma=\exp(T)$.
  Hence we have 
  \begin{equation}\label{app:3}
     \tilde\phi(C_{\gamma})=\Hyp^3\cap Z_{\gamma}=
           \exp\left(i\R T\right),
  \end{equation}
  because $\exp\left(i\Su(2)\right)=\Hyp^3$.

  Suppose now that $\Gamma$ is not abelian.
  Then there exist $\gamma$, $\gamma'\in \Gamma$
  such that $\gamma\!\cdot\! \gamma'\ne \gamma'\!\cdot \!\gamma$.
  Set $\gamma=\exp(T)$ and $\gamma'=\exp(T')$,
  where $T$, $T'\in \Su(2)$.
  Then we have $i {\R}T\cap i {\R}T'=\{0\}$.
  It is well-known that the restriction of the exponential map
  $\exp|_{i\Su(2)}\colon{}i\Su(2)\to\Hyp^3$ is bijective.
  Hence we have
  \[
      \tilde\phi(C_\gamma)\cap\tilde\phi( C_{\gamma'})=
      \exp\left(i{\R}T\right)\cap 
      \exp\left(i{\R}T'\right)=\{\id\}.
  \]
  By \eqref{app:1}, \eqref{app:2} and \eqref{app:3},
  we have
  \[
      \phi(I_{\Gamma})=\{\id\} \qquad 
      \mbox{(if $\Gamma$ is not abelian)}.
  \]
  Next we consider the case $\Gamma$ is abelian.
  If $\Gamma\subset \{\pm\id\}$, then obviously
  \[
      \phi(I_{\Gamma})=\Hyp^3.
  \]
  Suppose $\Gamma\not\subset \{\pm\id\}$. Then 
  there exists $\gamma\in \Gamma$ such that $\gamma\ne \pm \id$.
  We set $\gamma=\exp{T}$\,\,($T\in \Su(2)$).
  Since $\exp(\R T)$ is a maximal abelian subgroup
  containing $\gamma$, we have $\Gamma \subset \exp(\R T)$.
  Then by \eqref{app:3}, we have
  \[
      \phi(I_{\Gamma})=\exp(i\R T).
  \]
\end{Proof}

\par
\setlength{\parindent}{0in}
\footnotesize
\vspace{3mm}
{\sc Masaaki Umehara}:\,Department of Mathematics,
Graduate School of Science,
Osaka University,
Toyonaka, Osaka 560,
Japan.\,
{\it E-mail\/}:
umehara@math.wani.osaka-u.ac.jp
\par
\vspace{3mm}
{\sc Kotaro Yamada}:\,
Department of Mathematics,
Faculty of Science,
Kumamoto University,
Kumamoto 860,
Japan.\,
{\it E-mail\/}: 
kotaro@gpo.kumamoto-u.ac.jp

\end{document}